\documentstyle[12pt,a4paper,amssymb]{amsart}
\theoremstyle{plain}
\newtheorem{theorem}{Theorem}[section]
\newtheorem{corollary}[theorem]{Corollary}
\newtheorem{lemma}[theorem]{Lemma}

\theoremstyle{definition}

\theoremstyle{remark}

\theoremstyle{plain}

\theoremstyle{definition}

\theoremstyle{remark}

\numberwithin{equation}{section}

\newcommand{\mcal}[1] {{\mathcal{#1}}}
\newcommand{\mbb}[1] {{\mathbb{#1}}}
\newcommand{\bs}[1] {{\boldsymbol{#1}}}

\newcommand{\ess}{\sigma_{\mbox{\footnotesize ess}}}

\newcommand{\sing}{\sigma_{\mbox{\footnotesize sc}}}

\def \lnorm#1\rnorm {\vphantom{#1}\left\|\smash{#1}\right\|}
\def \lmod#1\rmod {\vphantom{#1}\left|\smash{#1}\right|}

\renewcommand \phi {\varphi}
\renewcommand \emptyset {\varnothing}

\newcommand \bp {{\bf p}}
\newcommand \bb {{\bf b}}
\newcommand \bx {{\bf x}}
\newcommand \be {{\bf e}}
\newcommand \bn {{\bf n}}
\newcommand \bB {{\bf B}}
\newcommand \bk {{\bf k}}
\newcommand \bK {{\bf K}}
\newcommand \bH {{\bf H}}
\newcommand \bJ {{\bf J}}

\begin{document}


\title[On the virial theorem for Brown-Ravenhall operator]%
   {On the virial theorem for the relativistic operator of
   Brown and Ravenhall, and the absence of embedded eigenvalues}

\date{\today}

\author{A.~ A.~Balinsky}
\address{School of Mathematics\\
   University of Wales, Cardiff\\
   23 Senghennydd Road, P O Box 926\\
   Cardiff CF2 4YH, Great Britain}
\email[A.~ A.~Balinsky]{BalinskyA@@cardiff.ac.uk}
\author{W.~D.~Evans}
\email[W.~D.~Evans]{EvansWD@@cardiff.ac.uk}

\keywords{Brown-Ravenhall operator, virial theorem, essential
spectra, eigenvalues}
\subjclass{Primary: 47G, 81Q10; Secondary: 47N50}

\begin{abstract}
A virial theorem is established for the operator
proposed by Brown and Ravenhall as a model for
relativistic one-electron atoms. As a consequence,
it is proved that the operator has no eigenvalues
greater than $\max(m c^2, 2 \alpha Z - \frac{1}{2})$,
where $\alpha$ is the fine structure constant,
for all values of the nuclear charge $Z$ below
the critical value $Z_c$: in particular there are no
eigenvalues embedded in the essential spectrum when
$Z \leq 3/4 \alpha$. Implications for the operators
in the partial wave decomposition are also described.
\end{abstract}

\maketitle


\section{Introduction}

\

	The formal operator proposed by Brown and Ravenhall in
\cite{BR}  to include relativistic effects in the description of
an electron in the field of its nucleus is of the form
\begin{equation}\label{oper}
{\bB} := \Lambda_+ \biggl(D_0 - \frac{e^2 Z}{\lmod \cdotp \rmod}\biggr)
\Lambda_+ .
\end{equation}
%
In (\ref{oper}), the notation is as follows (see \cite{EPS}):
\begin{itemize}
\item $D_0$ is the free Dirac operator
\begin{displaymath}
D_0= c {\bs{\alpha}} \cdotp \frac{\hslash}{i} \ \nabla  + m
c^2 \beta \equiv \sum_{j=1}^3 c \frac{\hslash}{i} \alpha_j
\frac{\partial}{\partial x_j} + m c^2 \beta , \end{displaymath}
where $\bs{\alpha}=(\alpha_1,\alpha_2,\alpha_3)$ and $\beta$ are the
Dirac matrices given by
\begin{displaymath}
\alpha_j=
\begin{pmatrix}
0 & \sigma_j\\
\sigma_j & 0
\end{pmatrix} , \ \
\beta=
\begin{pmatrix}
1_2 & 0_2\\
0_2 & -1_2
\end{pmatrix}
\end{displaymath}
with $\bs{0_2} , \bs{1_2}$ the zero and unit $2 \times 2$ matrices
respectively and $\sigma_j$ the Pauli matrices
\begin{displaymath}
\sigma_1=
\begin{pmatrix}
0 & 1\\
1 & 0
\end{pmatrix}, \ \
\sigma_2=
\begin{pmatrix}
0 & -i\\
i & 0
\end{pmatrix}, \ \
\sigma_3=
\begin{pmatrix}
1 & 0\\
0 & -1
\end{pmatrix}.
\end{displaymath}
\item $\Lambda_+$ denotes the projection of $\operatorname{L}^2
({\mathbb{R}}^3)\otimes {\mathbb{C}}^4$ onto the positive
spectral subspace of $D_0$, that is $\chi_{(0,\infty)} (D_0)$,
where $\chi_{(0,\infty)}$ is the characteristic function of
$(0,\infty)$. If we set
\begin{displaymath}
\widehat{f}(\bp) \equiv \mcal{F}(f) (\bp) = \biggl( \frac{1}{2
\pi \hslash}\biggr)^{3/2} \int_{{\mbb{R}}^3} e^{-i \ \bx \cdotp \bp /
\hslash} f(\bx) \operatorname{d} \bx
\end{displaymath}
for the Fourier transform of $f$, then it follows that
\begin{displaymath}
(\Lambda_+ f)^{\wedge} (\bp) = \Lambda_+ (\bp) \widehat{f} (\bp),
\end{displaymath}
where
\begin{equation}\label{proj}
\Lambda_+ (\bp)= \frac{1}{2} + \frac{c \bs{\alpha} \cdotp \bp
+ m c^2 \beta}{2 \be (p)}, \ \  \be (p) = \sqrt{c^2 p^2 + m^2
c^4}
\end{equation}
with $p=\lmod \bp\rmod$.
%
\item $2 \pi \hslash$ is Planck's constant, $c$ the velocity of
light, $m$ the electron mass, $-e$ the electron charge, and $Z$
the nuclear charge.
\end{itemize}

	The underlying Hilbert space in which ${\bB}$ acts is
\begin{equation}\label{space}
\mcal{H} = \Lambda_+(\operatorname{L}^2
({\mathbb{R}}^3)\otimes {\mathbb{C}}^4).
\end{equation}
The Fourier transform of any spinor $\psi$ in the positive
spectral subspace of $D_0$ can be written
\begin{equation}\label{positive}
\widehat{\psi}(\bp)= \frac{1}{\bn(p)}
\begin{pmatrix}
[ \be (p) + \be (0)] \ u (\bp)\\
c (\bp \cdotp \bs{\sigma} ) \ u (\bp)
\end{pmatrix},
\end{equation}
where $u \in \operatorname{L}^2
({\mathbb{R}}^3)\otimes {\mathbb{C}}^2$, a Pauli spinor, and
$\bn (p) = [ 2 \be (p)(\be (p) + \be (0))]^{1/2}$.   Conversely
any Dirac spinor  of the form (\ref{positive})\ is in the
image of $\mcal{H}$ under the Fourier transform. We than have
formally that, if $v$ is the Pauli spinor related by
(\ref{positive}) to  a Dirac spinor $\phi$,
\begin{equation}\label{form1} (\phi, \bB
\psi ) = \bs{\beta} [v,u] : = \be [v,u] - \gamma \bk [v,u],
\end{equation}
where  $(\cdotp,\cdotp)$ is the inner-product on
$ (\operatorname{L}^2 ({\mathbb{R}}^3)\otimes {\mathbb{C}}^4)$
(being linear in the second argument),
\begin{equation}\label{e-form}
\be [v,u] = \int_{{\mbb{R}}^3} \be  (p) \ v (\bp)^{\ast} \  u
(\bp) \operatorname{d} \bp ,
\end{equation}
\begin{equation}\label{k-form}
\bk [v,u] = \iint\limits_{{\mbb{R}}^3 \times {\mbb{R}}^3} v(\bp
')^\ast \bK (\bp',\bp) u(\bp)\operatorname{d} \bp'
\operatorname{d} \bp,
\end{equation}
and
\begin{equation}\label{gamma}
\gamma = \frac{\alpha c Z}{2 \pi^2}.
\end{equation}
In (\ref{e-form})-(\ref{gamma}), $\ast$ denotes the Hermitian
conjugate, $\alpha=e^2/(\hslash c)$ is Sommerfeld's fine
structure constant and the kernel $\bK$ in (\ref{k-form}) is the
$2 \times 2$ matrix-valued function
\begin{equation}\label{kernel}
\bK (\bp',\bp)= \frac{[\be (p')+\be(0)][\be(p)+\be(0)]
{\large\bs{1_2}} + c^2
(\bp' \cdotp \bs{\sigma})(\bp \cdotp \bs{\sigma})}{\bn (p') \lmod
\bp - \bp'\rmod^2 \bn (p)}.
\end{equation}

	The form $\be [\cdotp]$ with domain $\operatorname{L}^2
({\mathbb{R}}^3; \sqrt{1+ p^2} \operatorname{d} \bp)\otimes
{\mathbb{C}}^2$ is closed and non-negative in
$\operatorname{L}^2
({\mathbb{R}}^3)\otimes
{\mathbb{C}}^2$.
In \cite{BE} and \cite{TixPos} it was proved that for
$u \in \operatorname{L}^2
({\mathbb{R}}^3; \sqrt{1+ p^2} \operatorname{d} \bp)\otimes
{\mathbb{C}}^2$ and $ Z \leq Z_c = 2/[(\frac{\pi}{2} +
\frac{2}{\pi}) \alpha]$,
\begin{equation}\label{ineq1}
\gamma \bk[u,u] \leq \frac{Z}{Z_c} \be [u,u].
\end{equation}
It had  earlier been established in \cite{EPS} that
$\bs{\beta}$ is bounded below if and only if $Z \leq Z_c$,
confirming a prediction of  Hardenkopf and Sucher
\cite{HS} based on numerical considerations. The strict
positivity of $\bs{\beta}$ is proved in both \cite{BE}
\cite{TixPos}, a positive lower bound being exhibited in
\cite{TixPos} even for $Z=Z_c$. If $Z < Z_c$, it follows from
(\ref{form1}) and (\ref{ineq1}) that $\bs{\beta}$ is a closed
positive form with domain $\mcal{D}(\bs{\beta})=
\operatorname{L}^2 ({\mathbb{R}}^3; \sqrt{1+ p^2}
\operatorname{d} \bp)\otimes {\mathbb{C}}^2$.There is therefore
defined a positive self-adjoint operator $\bb$ in
$\operatorname{L}^2 ({\mathbb{R}}^3)\otimes {\mathbb{C}}^2 $
satisfying
\begin{equation}\label{b-oper} (v, \bb u) =
\bs{\beta} [v,u], \ \ u \in \mcal{D}(\bb), \ v \in
\mcal{D}(\bs{\beta}),
\end{equation}
where the domain  $\mcal{D}(\bb)$ of $\bb$ is dense in its
``form domain'' $\mcal{D}(\bs{\beta})$. We shall follow
\cite{EPS} and refer to this operator $\bb$ as the
Brown-Ravenhall operator, it being associated to the original
operator $\bB$ via (\ref{form1}).

	If $Z=Z_c$, the operator $\bb$ is defined as the Friedrichs
operator associated with the closure of the form $\bs{\beta}$
restricted to rapidly decreasing Pauli spinors. The description
of the operator $\bb$ in this case is studied in \cite{TixSelf}. We
assume in this paper that $Z < Z_c$, so that $\bb$ is defined in
(\ref{b-oper}). In this case, it follows from
\cite[Theorem~2]{EPS} that $\ess(\bb)=[m c^2,\infty)$ and
$\sing (\bb)=\emptyset$, where $\ess$ and $\sing$ are
respectively the essential and singular continuous spectra.
We shall establish a virial theorem for $\bb$ which will imply
that for $Z \leq Z_c'= \frac{3}{4 \alpha}$
there are no eigenvalues embedded in $[m c^2,\infty)$, and
hence the spectrum in $[m c^2,\infty)$ is absolutely continuous.

The number $Z_c'$ has a spectral significance in \cite{TixSelf}
concerning the self-adjointness of operators $b_{l,s}$
in the partial wave decomposition of $\bB$ and $\bb$:
in $b_{l,s}$, $l$ denotes the angular momentum chanel and $s$
the spin, and in \cite{TixSelf} the operators are given the domain
$\operatorname{L}^2 (0,\infty; [\sqrt{1+ p^2}]
\operatorname{d} p)$.
Tix proves that for $(l,s) \neq (0,1/2)$ or $(1,-1/2)$, the
operators $b_{l,s}$ are all self-adjoint for $Z < Z_c$, but
$b_{0,1/2}$ and $b_{1,-1/2}$ are self-adjoint if $Z < Z_c'$,
essentially self-adjoint when $Z = Z_c'$ and are symmetric with
 a one-parameter family of self-adjoint extensions when $ Z_c' <
 Z \leq Z_c$. We also prove that for all $(l,s) \neq (0, 1/2) $
 or $(1, -1/2)$, the $b_{l,s}$ have no eigenvalues in $[mc^2,
\infty )$ for $Z < Z_c$, but only for $Z \leq Z_c'$ in the case
$(l,s) = (1, - 1/2)$. The operator $b_{0,1/2}$ (its Friedrichs
extension when $Z \geq Z_c')$  has no eigenvalues in $[mc^2,
\infty )$ for the whole range $Z < Z_c$.

\

\section{The Virial Theorem}

\

	We begin by proving an abstract virial theorem.
\begin{lemma}\label{virial-abs}
Let $U(a),  \ a \in \mbb{R}_+$, be a one parameter family of
unitary operators on a Hilbert space $\mcal{H}$ which converges
strongly to the identity as $a \to 1$. Let $T$ be a self-adjoint
operator in $\mcal{H}$ and $T_a  =  f(a) \ U(a) T U(a)^{-1}$,
where $f(1)=1$ and $f'(1)$ exists. If $\phi \in \mcal{D}(T)\cap
\mcal{D} (T_a)$ is an eigenvector of $T$ corresponding to an
eigenvalue $\lambda$ then
\begin{equation}\label{lim-virial}
\lim\limits_{a \to 1} \biggl( \phi_a, \biggl[
\frac{T_a-T}{a-1}\biggr] \phi \biggr)=\lambda f'(1)\lnorm
\phi\rnorm^2,
\end{equation}
where $\phi_a=U(a)\phi$.
\end{lemma}

\

\begin{pf}
From $T \phi = \lambda \phi$ we have $T_a \phi_a = \lambda f(a)
\phi_a$. Hence
\begin{displaymath}
(\phi_a, T \phi) = \lambda (\phi_a,\phi)
\end{displaymath}
and
\begin{displaymath}
(T_a \phi_a,\phi) = \lambda f(a) (\phi_a, \phi).
\end{displaymath}
Consequently
\begin{displaymath}
\biggl( \phi_a, \biggl[
\frac{T_a-T}{a-1}\biggr] \phi \biggr)=\lambda \biggl[
\frac{f(a)-1}{a-1}\biggr] (\phi_a,\phi)
\end{displaymath}
and the result follows as allowing $a \to 1$.
\end{pf}

	Our main results are the following Theorem and Corollary
	which result from the application of Lemma~\ref{virial-abs}
	to the Brown-Ravenhall operator $\bb$.

\begin{theorem}\label{main-virial}
Let $Z < Z_c$ and let $\bb$ be defined by (\ref{b-oper})
 in $\operatorname{L}^2 ({\mathbb{R}}^3)\otimes
{\mathbb{C}}^2$  with form domain
$\mcal{D}(\bs{\beta}) = \operatorname{L}^2  ({\mathbb{R}}^3;
\sqrt{1+ p^2} \operatorname{d} \bp)\otimes {\mathbb{C}}^2$.
If $\lambda$ is an eigenvalue of $\bb$ with eigenvector $\phi$,
then
\begin{multline}
\lambda \lnorm \phi\rnorm^2 = \int_{\mbb{R}^3}
\frac{\be(0)^2}{\be (p)} \lmod \phi(\bp)\rmod^2 \operatorname{d}
\bp \\
 -\frac{\gamma}{2} \be(0) \iint\limits_{{\mbb{R}}^3 \times
{\mbb{R}}^3} \phi^\ast(\bp ') \bK^1 (\bp',\bp)
\phi(\bp) \biggl[
\frac{1}{\be(p)}-\frac{\be(0)}{\be(p)^2}+
\frac{1}{\be(p')}-\frac{\be(0)}{\be(p')^2}
\biggr]\operatorname{d} \bp' \operatorname{d} \bp \\
+
\frac{\gamma}{2} \be(0) \iint\limits_{{\mbb{R}}^3 \times
{\mbb{R}}^3} \phi^\ast(\bp ') \bK^2 (\bp',\bp)
\phi(\bp) \biggl[
\frac{1}{\be(p)}+\frac{\be(0)}{\be(p)^2}+
\frac{1}{\be(p')}+\frac{\be(0)}{\be(p')^2}
\biggr] \operatorname{d} \bp' \operatorname{d} \bp ,
\label{virial_BS}
\end{multline}
where
\begin{equation}\label{k1}
\bK^1 (\bp',\bp)= \frac{[\be (p')+\be(0)][\be(p)+\be(0)]
{\large\bs{1_2}}}{\bn (p') \lmod \bp - \bp'\rmod^2 \bn (p)},
\end{equation}
\begin{equation}\label{k2}
\bK^2 (\bp',\bp)= \frac{c^2 (\bp' \cdotp \bs{\sigma})(\bp \cdotp
\bs{\sigma})}{\bn (p') \lmod \bp - \bp'\rmod^2 \bn (p)}.
\end{equation}
\end{theorem}

\

\begin{corollary}\label{main-corol}
Let $Z < Z_c$ and $\bb \phi = \lambda \phi$. Then
\begin{multline}
\biggl( \frac{\lambda}{\be(0)}-1 \biggr)  \int_{\mbb{R}^3}
 \lmod \phi(\bp)\rmod^2  \biggl\{ 1-
 \frac{\be(0)}{\be(p)} +  \frac{\be(0)^2}{\be(p)^2}
 \biggr\}\operatorname{d} \bp \\
 =
  \gamma
 \iint\limits_{{\mbb{R}}^3 \times {\mbb{R}}^3} \phi^\ast(\bp ')
 \bK^2 (\bp',\bp) \phi(\bp) \biggl[
 \frac{1}{\be(p')} +
\frac{1}{\be(p)}
\biggr]\operatorname{d} \bp' \operatorname{d} \bp \\
-
\int_{\mbb{R}^3}
 \lmod \phi(\bp)\rmod^2
 \frac{[\be(p) - \be(0)] [2  \be(p) - \be(0)]}{\be(p)^2}
 \operatorname{d} \bp.
\label{corol_virial_BS}
\end{multline}
\end{corollary}

\

	Before proving Theorem~\ref{main-virial} and
	Corollary~\ref{main-corol}, we  need the following lemma.
	We shall hereafter in this section
	write $\bb_m, \ \be_m, \  \bK_m,$   to indicate the
	dependence on $m$.

\begin{lemma}\label{bounded}
For any $m \in \mbb{R}_+$, $\bb_m$ and $\bb_0$
have the same domain and
 $\bb_m - \bb_0$ can be extended to a
bounded operator on $\operatorname{L}^2  ({\mathbb{R}}^3)
\otimes {\mathbb{C}}^2$.
\end{lemma}

\

\begin{pf}
The result is essentially proved in  \cite[Theorem~1]{TixSelf}
but we give a short direct proof for completeness. The
following estimates are readily verified ( cf \cite{TixSelf}):
\begin{displaymath}
0 \leq \be_m (p) - \be_0 (p) \leq m c^2,
\end{displaymath}
\begin{multline*}
\lmod \bK_m (\bp',\bp) - \bK_0 (\bp',\bp)\rmod
\leq \frac{1}{\lmod \bp - \bp'\rmod^2}
\biggl\{ \frac{m}{2 \be_m (p)} + \frac{m}{2 \be_m (p')}  +
\frac{m^2}{4 \be_m (p)  \be_m (p')}\biggr\}  \\
\leq \frac{k(m)}{\lmod \bp - \bp'\rmod^2} \biggl\{
\frac{1}{\be_m (p)} + \frac{1}{\be_m (p')}\biggr\},
\end{multline*}
where $k(m)$ is a constant depending on $m$.
Hence, for $\phi \in
\mcal{D}(\bs{\beta}_m) =  \mcal{D}(\bs{\beta}_0)$,
\begin{multline*}
\lmod \bs{\beta}_m [\phi ]  -  \bs{\beta}_0 [\phi ] \rmod
\leq m c^2 \lnorm
\phi\rnorm^2 \\ + \gamma k(m) \iint\limits_{{\mbb{R}}^3 \times
{\mbb{R}}^3} \frac{1}{\lmod \bp - \bp'\rmod^2} \biggl\{
\frac{1}{\be_m (p)} + \frac{1}{\be_m (p')}\biggr\} \lmod
\phi(\bp)\rmod \lmod \phi(\bp')\rmod \operatorname{d} \bp'
\operatorname{d} \bp . \
\end{multline*}
By the
Cauchy-Schwarz inequality, the last integral is no greater than
\begin{displaymath}
\iint\limits_{{\mbb{R}}^3 \times
{\mbb{R}}^3} \lmod \phi(\bp)\rmod^2
 \frac{1}{\lmod \bp - \bp'\rmod^2}
\biggl\{
\frac{1}{\be_m (p)} + \frac{1}{\be_m (p')}\biggr\}
\frac{h(\bp')}{h(\bp)} \operatorname{d} \bp' \operatorname{d}
\bp ,
\end{displaymath}
where we choose $h(\bp)= \lmod \bp\rmod^{-\alpha}, \ 1 < \alpha
<2$. Then (see \cite[p124]{LL})
\begin{multline*}
\int_{{\mbb{R}}^3}
 \frac{1}{\lmod \bp - \bp'\rmod^2}
\frac{h(\bp')}{h(\bp) \be_m(p)} \operatorname{d} \bp'
= \frac{\lmod \bp\rmod^\alpha}{\be_m(p)}
\int_{{\mbb{R}}^3}
 \frac{1}{\lmod \bp - \bp'\rmod^2}
\frac{1}{\lmod \bp'\rmod^\alpha} \operatorname{d} \bp' \\
\leq \frac{\lmod \bp\rmod^{\alpha-1}}{c}
\int_{{\mbb{R}}^3}
 \frac{1}{\lmod \bp - \bp'\rmod^2}
\frac{1}{\lmod \bp'\rmod^\alpha} \operatorname{d} \bp' =O(1) \ \
\mbox{if} \ \ 1< \alpha <3,
\end{multline*}
\begin{multline*}
\int_{{\mbb{R}}^3}
 \frac{1}{\lmod \bp - \bp'\rmod^2}
\frac{h(\bp')}{h(\bp) \be_m(p')} \operatorname{d} \bp'
\leq \frac{\lmod \bp\rmod^\alpha}{c}
\int_{{\mbb{R}}^3}
 \frac{1}{\lmod \bp - \bp'\rmod^2}
\frac{1}{\lmod \bp'\rmod^{\alpha+1}} \operatorname{d} \bp' =O(1) \ \
\mbox{if} \ \ 0< \alpha <2.
\end{multline*}
Hence,
\begin{displaymath}
\lmod \bs{\beta}_m [\phi ]  -  \bs{\beta}_0 [\phi ] \rmod
\leq
k(m) \lnorm \phi\rnorm^2
\end{displaymath}
and $(\bs{\beta}_m   -  \bs{\beta}_0 ) [\cdotp , \cdotp]$
can be extended to a bounded sesquilinear form  on
$ \biggl\{  \operatorname{L}^2({\mathbb{R}}^3)\otimes
{\mathbb{C}}^2\biggr\} \times
\biggl\{  \operatorname{L}^2({\mathbb{R}}^3)\otimes
{\mathbb{C}}^2\biggr\}$ .
It follows that $\bb_m$ and $\bb_0$ have the same domain and
this in turn implies the rest of the lemma.
\end{pf}

\

\begin{pf*}{Proof of Theorem~\ref{main-virial}.}
We  apply Lemma~\ref{virial-abs} to $T=\bb_m$,
$f(a)=a$ and $U(a)$ defined by $U(a) \phi(\bp) =
a^{-3/2} \phi(\bp /a) =: \phi_a (\bp)$. Then
$U(a)\to \operatorname{I}$ strongly in
$\operatorname{L}^2({\mathbb{R}}^3)\otimes {\mathbb{C}}^2$ and
\begin{displaymath}
T_a := a U(a) \bb_m U(a)^{-1} = \bb_{m a}.
\end{displaymath}
By Lemma~\ref{bounded}, an eigenvector $\phi$ of $\bb_m$, lies
in $\mcal{D}(\bb_{m a})$ for any $a$. Thus, it remains to
evaluate the limit on the left-hand side of (\ref{lim-virial}).
We have
\begin{multline}
\biggl( \phi_a , \biggl[ \frac{\bb_{ma} - \bb_m}{a
- 1}\biggr] \phi \biggr)= \int_{\mbb{R}^3} \biggl[
\frac{\be_{ma}(p)-\be_m (p)}{a - 1}\biggr] \lmod
\phi(\bp)\rmod^2 \operatorname{d} \bp \\
- \gamma  \iint\limits_{{\mbb{R}}^3 \times {\mbb{R}}^3}
\phi_a(\bp ')^\ast \biggl[ \frac{\bK_{m a}
(\bp',\bp) - \bK_m (\bp',\bp)}{a-1}\biggr] \phi(\bp)\operatorname{d} \bp'
\operatorname{d} \bp =: I_1 - \gamma I_2 ,
\label{virial_proof}
\end{multline}
and (\ref{virial_BS})
will follow if we can justify taking the limit as $a \to 1$
under the integral signs on the right-hand side of
(\ref{virial_proof}).

  In $I_1$ the integrand is majorised by $\lmod \phi
  (\bp)\rmod^2$ and so the dominated convergence theorem applies.
  We write $I_2$ as

\begin{displaymath}
 \iint\limits_{{\mbb{R}}^3 \times {\mbb{R}}^3}
[\phi_a(\bp ')- \phi(\bp ')]^\ast \biggl[ \frac{\bK_{m a}
(\bp',\bp) - \bK_m (\bp',\bp)}{a-1}\biggr] \phi(\bp)\operatorname{d} \bp'
\operatorname{d} \bp
\end{displaymath}
 \begin{displaymath}
+ \iint\limits_{{\mbb{R}}^3 \times {\mbb{R}}^3}
\phi(\bp ')^\ast \biggl[ \frac{\bK_{m a}
(\bp',\bp) - \bK_m (\bp',\bp)}{a-1}\biggr] \phi(\bp)\operatorname{d} \bp'
\operatorname{d} \bp = : I_3 + I_4,
 \end{displaymath}
and use the readily verified estimate
\begin{displaymath}
\biggl|  \frac{\bK_{m a}
(\bp',\bp) - \bK_m (\bp',\bp)}{a-1}\biggr| \leq k_a/\lmod
\bp - \bp'\rmod^2 ,
\end{displaymath}
where  $k_a \to 1$ as $a \to 1$. Thus $I_3$ bounded by
\begin{displaymath}
  \iint\limits_{{\mbb{R}}^3 \times {\mbb{R}}^3}
\lmod \phi_a (\bp') - \phi (\bp')\rmod
 \lmod \phi (\bp) \rmod
 \frac{\operatorname{d}
\bp' \operatorname{d} \bp}{\lmod \bp - \bp'\rmod^2}
\end{displaymath}
and, on using the Parseval identity and $\mcal{F} \bigl(
\frac{1}{\lmod \cdotp \rmod}\bigr)= \sqrt{\frac{2}{\pi}}
\frac{1}{\lmod \cdotp \rmod^2}$ (with $\mcal{F}$ now the
standard Fourier transform with $\hbar =1$),
this is equal to
\begin{displaymath}
2 \pi^2 \int_{{\mbb{R}}^3}
u(\bx) v(\bx)
\frac{\operatorname{d} \bx}{\lmod \bx\rmod},
\end{displaymath}
where $\widehat{u} (\bp)=
\lmod \phi_a (\bp) - \phi (\bp)\rmod$ and $\widehat{v}
(\bp)=  \lmod \phi (\bp)\rmod$,
\begin{displaymath}
\leq \pi^3  \biggl(
\int_{{\mbb{R}}^3} \lmod \phi_a (\bp) - \phi (\bp)\rmod^2 \lmod
\bp\rmod \operatorname{d} \bp \biggr)^{1/2} \biggl(
\int_{{\mbb{R}}^3} \lmod \phi (\bp)\rmod^2 \lmod \bp\rmod
\operatorname{d} \bp \biggr)^{1/2}
\end{displaymath}
by the Cauchy-Schwarz inequality and
Kato's inequality \cite[p.~307]{Kato}. Thus $I_3 \to 0$ as $a
\to 1$ since $\phi \in \operatorname{L}^2 ({\mathbb{R}}^3;
\sqrt{1+ p^2} \operatorname{d} \bp)\otimes {\mathbb{C}}^2$. In
$I_4$ the integrand is majorised by the function $\frac{1}{\lmod
\bp - \bp'\rmod^2} \lmod \phi(\bp')\rmod \lmod \phi(\bp)\rmod $
which is integrable by Kato's inequality and hence the dominated
convergence theorem applies.  \end{pf*}

\

\begin{pf*}{Proof of Corollary~\ref{main-corol}.} \ \
From  (\ref{virial_BS}) we have
\begin{multline}
\lambda \lnorm \phi\rnorm^2 = \be(0)^2 \int_{\mbb{R}^3}
 \lmod \phi(\bp)\rmod^2  \frac{\operatorname{d}
\bp  }{\be (p)}  \\
 - \gamma \be(0) \operatorname{Re} \biggl[
 \iint\limits_{{\mbb{R}}^3 \times
{\mbb{R}}^3} \phi^\ast(\bp ') \bK (\bp',\bp)
\phi(\bp) \biggl\{ \frac{1}{\be(p)}-\frac{\be(0)}{\be(p)^2} \biggr\}
\operatorname{d} \bp' \operatorname{d} \bp \biggr] \\
+
\gamma \be(0) \iint\limits_{{\mbb{R}}^3 \times
{\mbb{R}}^3} \phi^\ast(\bp ') \bK^2 (\bp',\bp)
\phi(\bp) \biggl\{
\frac{1}{\be(p)}+
\frac{1}{\be(p')} \biggr\}
\operatorname{d} \bp' \operatorname{d} \bp.
\label{corol_pf_1}
\end{multline}
%
Also, for all $\psi \in \mcal{D}(\bs{\beta})$,
\begin{multline}
 \int_{\mbb{R}^3}
 \be(p) \psi^{\ast}(\bp) \phi(\bp)  \operatorname{d}
\bp
 - \gamma
 \iint\limits_{{\mbb{R}}^3 \times
{\mbb{R}}^3} \psi^\ast(\bp ') \bK (\bp',\bp)
\phi(\bp)
\operatorname{d} \bp' \operatorname{d} \bp  \\
= \lambda \int_{\mbb{R}^3}
  \psi^{\ast}(\bp) \phi(\bp)  \operatorname{d}
\bp
\label{corol_pf_2}.
\end{multline}
We choose $\psi(\bp)= \biggl[
\frac{1}{\be(p)}-\frac{\be(0)}{\be(p)^2}\biggr] \phi (\bp)$
in (\ref{corol_pf_2}): clearly $\psi \in \mcal{D}(\bs{\beta})$.
Then
\begin{multline}\label{corol_pf_3}
\gamma
\iint\limits_{{\mbb{R}}^3 \times
{\mbb{R}}^3} \phi^\ast(\bp ') \bK (\bp',\bp)
\phi(\bp) \biggl\{ \frac{1}{\be(p)}-\frac{\be(0)}{\be(p)^2} \biggr\}
\operatorname{d} \bp' \operatorname{d} \bp \\  =
\int_{\mbb{R}^3}
[\be(p) - \lambda] \biggl\{
\frac{1}{\be(p)}-\frac{\be(0)}{\be(p)^2}\biggr\}
\lmod \phi (\bp)\rmod^2
\operatorname{d} \bp.
\end{multline}
On substituting (\ref{corol_pf_3}) in (\ref{corol_pf_1}),
\begin{multline*}
\lambda \lnorm \phi\rnorm^2 = \be(0)^2 \int_{\mbb{R}^3}
 \lmod \phi(\bp)\rmod^2  \frac{\operatorname{d}
\bp  }{\be (p)}
 - \be(0)   \int_{\mbb{R}^3}
[\be(p) - \lambda] \biggl\{
\frac{1}{\be(p)}-\frac{\be(0)}{\be(p)^2}\biggr\}
\lmod \phi (\bp)\rmod^2
\operatorname{d} \bp \\
+
\gamma \be(0) \iint\limits_{{\mbb{R}}^3 \times
{\mbb{R}}^3} \phi^\ast(\bp ') \bK^2 (\bp',\bp)
\phi(\bp) \biggl\{
\frac{1}{\be(p)}+
\frac{1}{\be(p')} \biggr\}
\operatorname{d} \bp' \operatorname{d} \bp,
\end{multline*}
whence
\begin{multline*}
 \lambda  \int_{\mbb{R}^3}
 \lmod \phi(\bp)\rmod^2  \biggl\{ 1-
 \frac{\be(0)}{\be(p)} +  \frac{\be(0)^2}{\be(p)^2}
 \biggr\}\operatorname{d} \bp =
 \be(0)  \int_{\mbb{R}^3}
 \lmod \phi(\bp)\rmod^2  \biggl\{
 \frac{\be(0)}{\be(p)} -1 +  \frac{\be(0)}{\be(p)}
 \biggr\}\operatorname{d} \bp \\
 +  \gamma \be(0)
 \iint\limits_{{\mbb{R}}^3 \times {\mbb{R}}^3} \phi^\ast(\bp ')
 \bK^2 (\bp',\bp) \phi(\bp) \biggl[
 \frac{1}{\be(p')} +
\frac{1}{\be(p)}
\biggr]\operatorname{d} \bp' \operatorname{d} \bp,
\end{multline*}
and
\begin{multline}
(\lambda - \be(0))  \int_{\mbb{R}^3}
 \lmod \phi(\bp)\rmod^2  \biggl\{ 1-
 \frac{\be(0)}{\be(p)} +  \frac{\be(0)^2}{\be(p)^2}
 \biggr\}\operatorname{d} \bp \\
 =
\be (0) \int_{\mbb{R}^3}
 \lmod \phi(\bp)\rmod^2
 \frac{[\be(p) - \be(0)] [ \be(0) - 2  \be(p) ]}{\be(p)^2}
 \operatorname{d} \bp \\
 +
  \gamma \be(0)
 \iint\limits_{{\mbb{R}}^3 \times {\mbb{R}}^3} \phi^\ast(\bp ')
 \bK^2 (\bp',\bp) \phi(\bp) \biggl[
 \frac{1}{\be(p')} +
\frac{1}{\be(p)}
\biggr]\operatorname{d} \bp' \operatorname{d} \bp .
\label{corol_pf_4}
\end{multline}
The corollary is therefore proved.
\end{pf*}

\

\section{The absence of embedded eigenvalues}

\

\begin{theorem}\label{th31}
If $Z < Z_c$, the operator $\bb$ has no eigenvalues in
$[\max \bigl\{ 1, 2 \alpha Z - 1/2 \bigr\} m c^2, \infty )$. In
particular, if $Z \leq Z_c' = 3/(4 \alpha)$, there are no
eigenvalues in $[m c^2, \infty )$ and the spectrum of $\bb$ is
absolutely continuous.
\end{theorem}

\

\begin{pf}
From (\ref{corol_virial_BS})
\begin{multline}
\biggl(\frac{\lambda}{\be(0)} - 1\biggr)  \int_{\mbb{R}^3}
 \lmod \phi(\bp)\rmod^2  \biggl\{ 1-
 \frac{\be(0)}{\be(p)} +  \frac{\be(0)^2}{\be(p)^2}
 \biggr\}\operatorname{d} \bp \\
 \leq
 \gamma
 \iint\limits_{{\mbb{R}}^3 \times {\mbb{R}}^3} \phi^\ast(\bp ')
 \bH (\bp',\bp)  \lmod \phi(\bp )\rmod \lmod \phi(\bp)\rmod
 \operatorname{d} \bp' \operatorname{d} \bp \\
 -
 \int_{\mbb{R}^3}
 \lmod \phi(\bp)\rmod^2
 \frac{[\be(p) - \be(0)] [  2 \be(p)- \be(0)  ]}{\be(p)^2}
 \operatorname{d} \bp ,
\label{abs31}
\end{multline}
where
\begin{displaymath}
  \bH (\bp',\bp)   =
\frac{c^2 p p'}{\bn(p') \bn(p) \lmod \bp - \bp'\rmod^2}
\biggl[ \frac{1}{\be(p')} +
\frac{1}{\be(p)}
\biggr].
\end{displaymath}
By the Cauchy-Schwarz inequality
\begin{equation}
 \iint\limits_{{\mbb{R}}^3 \times {\mbb{R}}^3}
  \bH (\bp',\bp)  \lmod \phi(\bp')\rmod \lmod
 \phi(\bp)\rmod \operatorname{d} \bp' \operatorname{d} \bp
\leq  c^2
\int_{\mbb{R}^3} \frac{p^2}{\bn^2(p)}
 \lmod \phi(\bp)\rmod^2
 \operatorname{d} \bp
\int_{\mbb{R}^3} \frac{h(\bp')}{h(\bp)}
 \bJ (\bp',\bp)
 \operatorname{d} \bp' ,
\label{abs32}
\end{equation}
where
\begin{displaymath}
   \bJ (\bp',\bp) =	   \frac{1}{\lmod \bp - \bp'\rmod^2} \biggl[
 \frac{1}{\be(p')} +
\frac{1}{\be(p)}
\biggr]
\end{displaymath}
and $h(\cdotp)$ is an arbitrary positive function. We make the
choice $h(\bp)= p^{-3/2}$ and use the following
(see \cite[p~124]{LL}):
$$
\int_{\mbb{R}^3}   \frac{1}{\lmod \bp - \bp'\rmod^2}
\frac{1}{\lmod \bp'\rmod^{3/2}}
\operatorname{d} \bp'   = 4 \pi^2 p^{-1/2},
$$
$$
\int_{\mbb{R}^3}   \frac{1}{\lmod \bp - \bp'\rmod^2}
\frac{1}{\lmod \bp'\rmod^{5/2}}
\operatorname{d} \bp'   = 4 \pi^2 p^{-3/2}.
$$
On substituting in (\ref{abs32}), we have
\begin{equation*}
 \iint\limits_{{\mbb{R}}^3 \times {\mbb{R}}^3}
  \bH (\bp',\bp)  \lmod \phi(\bp')\rmod \lmod
 \phi(\bp)\rmod \operatorname{d} \bp' \operatorname{d} \bp
< 4 \pi^2 c^2
\int_{\mbb{R}^3} \frac{p^2}{\bn^2(p)}
\biggl[ \frac{p}{\be(p)} + \frac{1}{c}\biggr]
\lmod \phi(\bp)\rmod^2
 \operatorname{d} \bp
\end{equation*}
and hence (\ref{abs31}) yields
\begin{multline*}
(\frac{\lambda}{\be(0)} - 1)  \int_{\mbb{R}^3}
 \lmod \phi(\bp)\rmod^2  \biggl\{ 1-
 \frac{\be(0)}{\be(p)} +  \frac{\be(0)^2}{\be(p)^2}
 \biggr\}\operatorname{d} \bp \\
 <
 \int_{\mbb{R}^3}
 \lmod \phi(\bp)\rmod^2
 \biggl\{
 4 \pi^2 c^2 \gamma
 \frac{p^2}{\bn^2(p)}
\biggl[ \frac{p}{\be(p)} + \frac{1}{c}\biggr] -
 \frac{[\be(p) - \be(0)] [  2 \be(p)- \be(0)
 ]}{\be(p)^2} \biggr\}
 \operatorname{d} \bp .
\end{multline*}
On replacing $\bp$ by $m c \bp$ and simplifying, we derive
\begin{equation}\label{abs33}
0 < \int_{\mbb{R}^3}
\lmod \phi(m c \bp)\rmod^2
\biggl(1 - \frac{1}{\sqrt{p^2 +1}} + \frac{1}{p^2 +1} \biggr)
\biggl(1 - \frac{\lambda}{m c ^2} + \Phi (p)\biggr)
\operatorname{d} \bp ,
\end{equation}
where
\begin{multline*}
\Phi (p)= \frac{p^2 (p + \sqrt{p^2 +1})}{(\sqrt{p^2 +1} +1)(p^2
+2 - \sqrt{p^2 +1})}
\biggl\{ 2 \pi^2 \biggl( \frac{\gamma}{c}\biggr) -
\frac{2 \sqrt{p^2 +1} - 1}{p + \sqrt{p^2 +1}}
\biggr\} \\
\leq
\frac{p^2 (p + \sqrt{p^2 +1})}{(\sqrt{p^2 +1} +1)(p^2
+2 - \sqrt{p^2 +1})}
\biggl\{ 2 \pi^2 \biggl( \frac{\gamma}{c}\biggr) -
\frac{3}{4}
\biggr\}
\end{multline*}
since $\min\limits_{\mbb{R}_+} \biggl\{
\frac{2 \sqrt{p^2 +1} - 1}{p + \sqrt{p^2 +1}}
\biggr\} = \frac{3}{4}$.  Thus, if $ 2 \pi^2 \bigl(
\frac{\gamma}{c}\bigr) \equiv \alpha Z \leq \frac{3}{4},$
(\ref{abs33}) implies that $\lambda < mc^2$. Since
$\sup\limits_{\mbb{R}_+} \biggl\{  \frac{p^2 (p + \sqrt{p^2
+1})}{(\sqrt{p^2 +1} +1)(p^2 +2 - \sqrt{p^2 +1})}\biggr\}=2$,
it follows from (\ref{abs33}) that, if
$ 2 \pi^2 \bigl(
\frac{\gamma}{c}\bigr) \equiv \alpha Z > \frac{3}{4},$
$\lambda < mc^2 [1+ 2 \alpha Z - \frac{3}{2}]$. The theorem is
therefore proved
\end{pf}

\

	Finally, we analyse the implications of the virial theorem
	for the operators $b_{l,s}$.

	In the partial wave decomposition of $\bB$ and $\bb$ (see
\cite{EPS}), spinors
$\phi \in\operatorname{L}^2~({\mathbb{R}}^3)\otimes
{\mathbb{C}}^2$ are expanded in terms of spherical spinors
$\Omega_{l,n,s}$
\begin{equation}\label{wave-decomp}
\phi (\bp) = \sum\limits_{(l,n,s) \in I} p^{-1} \phi_{l,n,s}(p)
\Omega_{l,n,s} (\omega)
\end{equation}
say, where $I= \{ (l,n,s) : l \in \mcal{N}_0 , n = -l -
\frac{1}{2}, \cdots, l + \frac{1}{2}, s = \pm \frac{1}{2} ,
\Omega_{l,n,s} \neq 0\}$ and
\begin{equation}\label{wave-norm}
\int_{\mbb{R}^3}
\lmod \phi (\bp)\rmod^2
\operatorname{d} \bp =
\sum\limits_{(l,n,s) \in I}
\int\limits_{0}^{\infty}
\lmod \phi_{l,n,s} (p)\rmod^2
\operatorname{d} p .
\end{equation}
On substituting in (\ref{form1}), we get
\begin{equation}\label{form_decomp_1}
\bs{\beta} [\phi , \psi ]=
\sum\limits_{(l,n,s) \in I}
\bs{\beta}_{l,s} [\phi_{l,n,s} , \psi_{l,n,s} ],
\end{equation}
where
\begin{multline}\label{form_decomp_2}
\bs{\beta}_{l,s} [\phi_{l,n,s} , \psi_{l,n,s} ]=
\int\limits_{0}^{\infty}
\be (p)\overline{\phi}_{l,n,s} (p) \psi_{l,n,s} (p)
\operatorname{d} p \\
- \frac{\alpha c Z}{\pi}
\int\limits_{0}^{\infty}\int\limits_{0}^{\infty}
\overline{\phi}_{l,n,s}(p') k_{l,s} (p',p) \psi_{l,n,s} (p)
\operatorname{d} p' \operatorname{d} p
\end{multline}
and, with $Q_l$ denoting Legendre functions of the second kind,
\begin{eqnarray}\label{wave-kernal}
k_{l,s} (p',p)& =&
\frac{[\be(p') +
\be(0)] Q_l \bigl( \frac{1}{2} \bigl[\frac{p'}{p}
+ \frac{p}{p'}\bigr]\bigr) [\be(p) +
\be(0)]}{\bn(p') \bn(p)}
+ \frac{c^2 p' Q_{l+ 2s} \bigl( \frac{1}{2}[\frac{p'}{p}
+ \frac{p}{p'}]\bigr) p }{\bn(p') \bn(p)}   \nonumber \\
& & \ \ \nonumber \\
& & =: k_{l,m}^1 (p',p) + k_{l,m}^2 (p',p)
\end{eqnarray}
say. For $Z < Z_c$, the forms with domain
$\operatorname{L}^2
(0, \infty; \sqrt{1+ p^2} \operatorname{d} p)$
are closed and positive in
$\operatorname{L}^2
(0, \infty)$, and we shall denote the associated self-adjoint
operators by $b_{l,s}$. These operators $b_{l,s}$ coincide with
the operators  $b_{l,s}$ of Tix in \cite{TixSelf} when the later
are self-adjoint,  but  are otherwise their Friedrichs
extensions. Note that it follows from \cite[Theorem2]{EPS} that
for all values of $l,s$, $\ess(b_{ls})= [mc^2,\infty)$.

\

\begin{theorem}\label{virial_for_decom}
For $(l,s) \neq (1, - \frac{1}{2}),$  the operators $b_{l,s}$
have no eigenvalues in $[m c ^2 , \infty)$ if $Z < Z_c$.  If
$(l,s) = (1, - \frac{1}{2}),$ $b_{l,s}$ has no eigenvalues in
$[\max \bigl\{ 1, 2 \alpha Z - 1/2 \bigr\} m c^2, \infty )$:
in particular $b_{1,-1/2}$ has no eigenvalues in
$[m c ^2 , \infty)$
if  $Z \leq Z_c' = \frac{3}{4 \alpha} $.
\end{theorem}

\begin{pf}
The analogue of (\ref{corol_virial_BS}) is
\begin{multline}
\biggl(\frac{\lambda}{\be(0)} - 1\biggr)  \int\limits_0^{\infty}
 \lmod \phi_{l,n,s}(p)\rmod^2  \biggl\{ 1-
 \frac{\be(0)}{\be(p)} +  \frac{\be(0)^2}{\be(p)^2}
 \biggr\}\operatorname{d} p \\
=  2 \pi
 \gamma
 \int\limits_0^{\infty} \int\limits_0^{\infty}
 \overline{\phi}_{l,n,s}(p ')
 \bk_{l,s}^2 (p',p)  \phi(p)
 \biggl[ \frac{1}{\be(p')} + \frac{1}{\be(p)}\biggr]
   \operatorname{d} p' \operatorname{d} p \\
-  \int\limits_0^{\infty}
 \lmod \phi_{l,n,s}(p)\rmod^2
 \frac{[\be(p) - \be(0)] [  2 \be(p)- \be(0)  ]}{\be(p)^2}
 \operatorname{d} p .
\label{abs39}
\end{multline}
Since $Q_0 (t) \geq \cdots \geq Q_l (t) \geq 0$ for all $t > 1$,
it follows that for all $(l,s)\neq (1, -1/2)$
\begin{displaymath}
 k^2_{l,s} (p',p) \leq \frac{c ^2 p' p}{\bn(p') \bn(p)} Q_1
 \biggl( \frac{1}{2}\biggl[ \frac{p'}{p} + \frac{p}{p'}\biggr]
 \biggr)
 \end{displaymath}
and hence, on writing $\phi$ for $\phi_{l,n,s}$ in (\ref{abs39}),
we get
\begin{multline}
\biggl(\frac{\lambda}{\be(0)} - 1\biggr)  \int\limits_0^{\infty}
 \lmod \phi(p)\rmod^2  \biggl\{ 1-
 \frac{\be(0)}{\be(p)} +  \frac{\be(0)^2}{\be(p)^2}
 \biggr\}\operatorname{d} p \\
\leq 2 \pi
 \gamma c^2
 \int\limits_0^{\infty}
 \lmod \phi(p)\rmod^2
 \frac{p^2}{\bn^2 (p) h(p)}
 \operatorname{d} p
 \int\limits_0^{\infty}
 h(p')
 Q_1
 \biggl( \frac{1}{2}\biggl[ \frac{p'}{p} + \frac{p}{p'}\biggr]
 \biggr)
\biggl( \frac{1}{\be(p')} + \frac{1}{\be(p)}\biggr)
   \operatorname{d} p'  \\
-  \int\limits_0^{\infty}
 \lmod \phi(p)\rmod^2
 \frac{[\be(p) - \be(0)] [  2 \be(p)- \be(0)  ]}{\be(p)^2}
 \operatorname{d} p ,
\label{abs310}
\end{multline}
for any positive function $h(\cdotp)$, on applying the
Cauchy-Schwarz inequality. We choose $h(t)=1/t$, so that in the
first term on the right-hand side of (\ref{abs310})we have
\begin{multline*}
 2 \pi
 \gamma c^2
 \int\limits_0^{\infty}
 \lmod \phi(p)\rmod^2
 \frac{p^2}{\bn^2 (p)}
 \biggl\{
  \frac{p}{\be(p)}
   \int\limits_0^{\infty}
 \frac{1}{p'}
 Q_1
 \biggl( \frac{1}{2}\biggl[ \frac{p'}{p} + \frac{p}{p'}\biggr]
 \biggr)
   \operatorname{d} p' \\
    +
   p \int\limits_0^{\infty}
 \frac{1}{p' \be(p')}
 Q_1
 \biggl( \frac{1}{2}\biggl[ \frac{p'}{p} + \frac{p}{p'}\biggr]
 \biggr)
   \operatorname{d} p'
   \biggr\}
 \operatorname{d} p \\
 =
 2 \pi
 \gamma c^2
 \int\limits_0^{\infty}
 \lmod \phi(p)\rmod^2
 \frac{p^2}{\bn^2 (p)}
 \biggl\{
I_1 + I_2
\biggr\} \operatorname{d} p
\end{multline*}
say. Let $g_1 (u)=  Q_1 \biggl(
 \frac{1}{2}\biggl[ u + \frac{1}{u}\biggr]
 \biggr) $. Then
 \begin{displaymath}
 c I_1 \leq \int\limits_{0}^{\infty} g_1 (u) \frac{
 \operatorname{d} u }{u} = 2,
\end{displaymath}
(see \cite[(3.8)]{BE} or \cite[\S~2.3]{EPS}), and
\begin{displaymath}
c I_2 \leq \int\limits_{0}^{\infty} g_1 (u) \frac{
 \operatorname{d} u }{u^2} =
 \int\limits_{0}^{\infty} g_1 (u)
 \operatorname{d} u  .
\end{displaymath}
We have (see \cite[(3.7)]{BE})
\begin{multline*}
 \int\limits_{0}^{1} g_1 (u)
 \operatorname{d} u =
 \frac{1}{2}  \int\limits_{0}^{1} u \ln \biggl( \frac{u
 + 1}{1 - u}\biggr) \operatorname{d} u
+ \frac{1}{2} \int\limits_{0}^{1} \frac{1}{u} \ln \biggl( \frac{u
 + 1}{1 - u}\biggr) \operatorname{d} u
 -
 \int\limits_{0}^{1}  \operatorname{d} u \\
 =\frac{1}{2} +  \frac{1}{2} \frac{\pi^2}{4} - 1 =
 \frac{\pi^2}{8} - \frac{1}{2}
\end{multline*}
  and
\begin{multline*}
 \int\limits_{1}^{\infty} g_1 (u)
 \operatorname{d} u =
 \frac{1}{2}  \int\limits_{1}^{\infty} \frac{1}{u} \ln \biggl(
 \frac{u + 1}{u-1}\biggr) \operatorname{d} u +
\int\limits_{1}^{\infty} \biggl\{\frac{u}{2} \ln \biggl( \frac{u
 + 1}{u-1}\biggr) - 1  \biggr\} \operatorname{d} u \\
 =\frac{\pi^2}{8} +  \biggl[ \frac{1}{4}(u^2 -1) \ln \biggl(
 \frac{u + 1}{ u-1}\biggr) - \frac{u}{2}\biggr]_1^{\infty}  =
 \frac{\pi^2}{8} +\frac{1}{2} .
 \end{multline*}
Hence $c I_2 \leq \frac{\pi^2}{4}$, and from (\ref{abs310})
\begin{multline*}
\biggl(\frac{\lambda}{\be(0)} - 1\biggr)  \int\limits_0^{\infty}
 \lmod \phi(p)\rmod^2  \biggl\{ 1-
 \frac{\be(0)}{\be(p)} +  \frac{\be(0)^2}{\be(p)^2}
 \biggr\}\operatorname{d} p \\
\leq 2 \pi
 \gamma c
 \int\limits_0^{\infty}
 \lmod \phi(p)\rmod^2
 \frac{p^2}{\bn^2 (p) } \bigl( \frac{\pi^2}{4} + 2\bigr)
\operatorname{d} p
-  \int\limits_0^{\infty}
 \lmod \phi(p)\rmod^2
 \frac{[\be(p) - \be(0)] [  2 \be(p)- \be(0)  ]}{\be(p)^2}
 \operatorname{d} p .
\end{multline*}
On replacing $p$ by $m c p$ and simplifying, we obtain
\begin{equation*}
0 < \int_{0}^{\infty}
\lmod \phi(m c \bp)\rmod^2
\biggl(1 - \frac{1}{\sqrt{p^2 +1}} + \frac{1}{p^2 +1} \biggr)
\biggl(1 - \frac{\lambda}{m c ^2} + \Psi (p)\biggr)
\operatorname{d} p ,
\end{equation*}
where
\begin{multline*}
\Psi (p)= \frac{p^2 (\sqrt{p^2 +1})}{(\sqrt{p^2 +1} +1)(p^2
+2 - \sqrt{p^2 +1})}
\biggl\{  \pi \biggl( \frac{\gamma}{c}\biggr)
\biggl(\frac{\pi^2}{4} + 2 \biggr) -
\frac{2 \sqrt{p^2 +1} - 1}{\sqrt{p^2 +1}}
\biggr\} \\
\leq
\frac{p^2 \sqrt{p^2 +1}}{(\sqrt{p^2 +1} +1)(p^2
+2 - \sqrt{p^2 +1})}
\biggl\{
\frac{2}{\pi^2 +4} \biggl( \frac{\pi^2}{4} + 2\biggr) -1
\biggr\} < 0
\end{multline*}
for $ 2 \pi^2 \frac{\gamma}{c} \equiv \alpha Z \leq \alpha Z_c$.
Hence $\lambda < m c ^2$.

	For the case $(l,s)=(1,-1/2)$, the Legendre function $Q_1$
in (\ref{abs310}) has to be replaced by $Q_0$, and we have to
consider
\begin{equation}\label{eigen_11}
 \frac{1}{h (p)}
\int\limits_{0}^{\infty} h(p') Q_0 \biggl(
 \frac{1}{2}\biggl[ \frac{p'}{p} + \frac{p}{p'}\biggr]
 \biggr)\biggr]
\biggl[
 \frac{1}{\be(p')} +
\frac{1}{\be(p)}
\biggr]       \operatorname{d} p' .
\end{equation}
We now make the choice $h(u)= 1/\sqrt{u}$, and, with
$g_0 (u)=  Q_0 \biggl(
 \frac{1}{2}\biggl[ u + \frac{1}{u}\biggr]
 \biggr) \equiv \ln \biggl| \frac{u + 1}{u - 1}\biggr|$,
(\ref{eigen_11}) becomes $J_1 + J_2$ say, where
 \begin{displaymath}
 J_1 = p \int\limits_{0}^{\infty}
 \frac{1}{\sqrt{u} \be(pu)} g_0 (u)     \operatorname{d} u
 < \frac{1}{c} \int\limits_{0}^{\infty}
   \frac{ g_0 (u)}{u^{3/2}}     \operatorname{d} u
 \end{displaymath}
 and
 \begin{displaymath}
 J_2 = \frac{p}{\be(p)} \int\limits_{0}^{\infty}
 \frac{1}{\sqrt{u}} g_0 (u)     \operatorname{d} u  .
 \end{displaymath}
  Since
 \begin{eqnarray*}
 \int\limits_{0}^{1}
   \frac{ g_0 (u)}{u^{1/2}} \operatorname{d} u=
   \int\limits_{1}^{\infty}
   \frac{ g_0 (u)}{u^{3/2}} \operatorname{d} u, \\
   \int\limits_{1}^{\infty}
   \frac{ g_0 (u)}{u^{1/2}} \operatorname{d} u=
   \int\limits_{0}^{1}
   \frac{ g_0 (u)}{u^{3/2}} \operatorname{d} u
 \end{eqnarray*}
 we have
 \begin{displaymath}
 \int\limits_{0}^{\infty}
   \frac{ g_0 (u)}{u^{1/2}} \operatorname{d} u=
   \int\limits_{0}^{\infty}
   \frac{ g_0 (u)}{u^{3/2}} \operatorname{d} u=
   \int\limits_{1}^{\infty}
    g_0 (u) \biggl( 1 + \frac{1}{u}\biggr)
	\frac{\operatorname{d} u}{\sqrt{u}}
\end{displaymath}
\begin{displaymath}
= 4 \int\limits_{1}^{\infty}
(\ln y) \frac{y}{(y^2 -1)^{3/2}}
\operatorname{d} y
\end{displaymath}
on setting $\frac{u + 1}{u - 1} = y$,
\begin{displaymath}
=  \int\limits_{1}^{\infty}
 (\ln z)  \frac{\operatorname{d} z }{(z-1)^{3/2}}
 =
2 \int\limits_{1}^{\infty}
   \frac{ 1}{z(z-1)^{1/2}} \operatorname{d} z = 2 \pi .
\end{displaymath}
Hence, on substituting in (\ref{abs310}) (with $Q_1$ replaced by
$Q_0$) we infer that
\begin{multline*}
\biggl(\frac{\lambda}{\be(0)} - 1\biggr)  \int\limits_0^{\infty}
 \lmod \phi(p)\rmod^2  \biggl\{ 1-
 \frac{\be(0)}{\be(p)} +  \frac{\be(0)^2}{\be(p)^2}
 \biggr\}\operatorname{d} p \\
<  2 \pi
 \gamma c
  \int\limits_0^{\infty}
	\lmod \phi(p)\rmod^2 \frac{p^2}{\bn^2 (p)}
	\biggl[ 2 \pi \biggl( \frac{c p}{\be (p)} +1\biggr)\biggr]
	 \operatorname{d} p
-  \int\limits_0^{\infty}
 \lmod \phi(p)\rmod^2
 \frac{[\be(p) - \be(0)] [  2 \be(p)- \be(0)  ]}{\be(p)^2}
 \operatorname{d} p .
\end{multline*}
On replacing $p$ by $m c p$ and simplifying, this yields
\begin{displaymath}
0 < \int_{0}^{\infty}
\lmod \phi(m c \bp)\rmod^2
\biggl(1 - \frac{1}{\sqrt{p^2 +1}} + \frac{1}{p^2 +1} \biggr)
\biggl(1 - \frac{\lambda}{m c ^2} + \Theta (p)\biggr)
\operatorname{d} p ,
\end{displaymath}
where
\begin{multline*}
\Theta (p)= \frac{p^2 ( p + \sqrt{p^2 +1})}{(\sqrt{p^2 +1}
+1)(p^2 +2 - \sqrt{p^2 +1})} \biggl\{ 2 \pi^2 \biggl(
\frac{\gamma}{c}\biggr) - \frac{2 \sqrt{p^2 +1} - 1}{p +
\sqrt{p^2 +1}} \biggr\} \\
\leq \frac{p^2 (\sqrt{p^2
+1} +p )}{(\sqrt{p^2 +1} +1)(p^2 +2 - \sqrt{p^2 +1})} \biggl\{
2 \pi^2 \biggl(
\frac{\gamma}{c}\biggr) - \frac{3}{4}
\biggr\} \leq 0
\end{multline*}
if $2 \pi^2 \bigl(
\frac{\gamma}{c}\bigr) = \alpha Z \leq \alpha Z_c'=3/4$.
Hence in
this case $\lambda < mc^2$. If $\alpha Z >3/4$, we have
\begin{displaymath}
\Theta (p) \leq 2 \bigl( \alpha Z - \frac{3}{4}\bigr)
\end{displaymath}
and $\lambda < mc^2 \bigl( \alpha Z - \frac{1}{2}\bigr)$.
The theorem is therefore proved.
\end{pf}

\

\section*{Acknowledgements}
The authors are grateful to the European Union for support
under the TMR grant FMRX-CT 96-0001.


\

\

\end{document}